\documentclass[leqno,11pt]{article}
\usepackage[spanish]{babel}
\usepackage[utf8]{inputenc}
\usepackage{amscd}
\usepackage{amsthm}
\usepackage{amssymb}
\usepackage{amsmath}
\usepackage{graphics}
\usepackage{graphicx}
\usepackage{verbatim}
\usepackage{color}
\newtheorem{mydef}{Definition}
\newtheorem{myexa}{Example}

\newtheorem{mytheo}{Theorem}
\newtheorem{mylemma}{Lemma}

\newtheorem{myprop}{Proposition}

\newtheorem*{myprem}{{\em Remarks}}
\newtheorem*{myproo}{{\em Proof}}

\newcommand{\pt}{\mbox{$\succ$\hspace{-1ex}$\longrightarrow$}}

\decimalpoint

\begin{document}

\begin{center}
	{\sc Consistency of the Bayes Estimator of a Regression Curve}\vspace{2ex}\\
	A.G. Nogales\vspace{2ex}\\
	Dpto. de Matem\'aticas, Universidad de Extremadura\\
	Avda. de Elvas, s/n, 06006--Badajoz, SPAIN.\\
	e-mail: nogales@unex.es
\end{center}
\vspace{.4cm}
\begin{quote}
	\hspace{\parindent} {\small {\sc Abstract.}  Strong consistency of the Bayes estimator of  a regression curve for the $L^1$-squared loss function is proved. It is also shown the convergence to 0 of the Bayes risk of this estimator both for the $L^1$ and $L^1$-squared loss functions. The Bayes estimator of a regression curve is the regression curve with respect to the posterior predictive distribution. 
	}
\end{quote}

\vspace{3ex}
\begin{itemize}
	\item[] \hspace*{-1cm} {\em AMS Subject Class.} (2020): 62F15, 62G08, 62G20.
	\item[] \hspace*{-1cm} {\em Key words and phrases: } Bayesian estimation of a regression curve, posterior predictive distribution.
\end{itemize}

\section{Introduction.}

Nogales (2022a) addresses the problem of estimation of a density from a Bayesian point of view and, under mild conditions, shows that the posterior predictive density is the Bayes estimator for the $L^1$-squared loss function, and Nogales (2022c) shows the strong consistency of this estimator. Nogales (2022b) deals, among others, with the problem of estimation of a regression curve and proves that the regression curve with respect to the posterior predicitive distribution is the Bayes estimator and, here, we wonder about its consistency. This is the main goal of the paper and the Theorem 1 below answers the question in the affirmative. 

The interested reader can find in  the papers mentioned above, and the references therein, more information on the problems of estimating a density or a regression curve, both from a frequentist and a Bayesian perspective, or about the usefulness of the posterior predictive distribution in Bayesian Inference and its calculation. We place special emphasis on the monographs Geisser (1993), Gelman et al. (2014) and Ghosal et al. (2017). 

Some examples are provided to illustrate the main result of the paper. For ease of reading we reproduce here an appendix from Nogales (2022a) to recall  the basic concepts of Bayesian inference but, mainly, to explain the (rather unusual) notation used in the paper.

\section{The framework.}

We recall from Nogales (2022b) the appropriate framework to address the problem, and update it to incorporate the required asymptotic flavor.

Let $(\Omega,\mathcal A,\{P_\theta\colon\theta\in (\Theta,\mathcal T,Q)\})$ be a Bayesian statistical experiment and $X_i:(\Omega,\mathcal A,\{P_\theta\colon\theta\in(\Theta,\mathcal T,Q)\})\rightarrow(\Omega_i,\mathcal A_i)$, $i=1,2$, two statistics. Consider the Bayesian experiment image of $(X_1,X_2)$
$$(\Omega_1\times\Omega_2,\mathcal A_1\times\mathcal  A_2,\{P_\theta^{(X_1,X_2)}\colon\theta\in(\Theta,\mathcal T,Q)\}).
$$
In the next, we shall suppose that $P^{(X_1,X_2)}(\theta,A_{12}):=P_\theta^{(X_1,X_2)}(A_{12})$, $\theta\in\Theta, \  A_{12}\in \mathcal A_1\times\mathcal  A_2$, is a Markov kernel. 

Let us write $R_\theta=P_\theta^{(X_1,X_2)}$ and $p_j(x):=x_j$ for $j=1,2$, $x:=(x_1,x_2)\in\Omega_1\times\Omega_2$. Hence
 \begin{gather*}
 	P_\theta^{X_1}=R_\theta^{p_1},\quad P_\theta^{X_2|X_1=x_1}=R_\theta^{p_2|p_1=x_1}\quad\hbox{and}\quad
 	E_{P_\theta}(X_2|X_1=x_1)=E_{R_\theta}(p_2|p_1=x_1).
 \end{gather*}

Given an integer $n$, for $m=n$ (resp. $m=\mathbb N$),  the Bayesian experiment corresponding to a $n$-sized sample (resp. an infinite sample) of the joint distribution of $(X_1,X_2)$ is 
$$\big((\Omega_1\times\Omega_2)^{m},(\mathcal A_1\times\mathcal  A_2)^{m},\big\{R_\theta^{m}\colon\theta\in(\Theta,\mathcal T,Q)\big\}\big)\qquad \hbox{(1)}
$$
We write $R^{m}(\theta,A'_{12,m})=R_\theta^{m}(A'_{12,m})$ for $A'_{12,m}\in (\mathcal A_1\times\mathcal  A_2)^{m}$ and 
$$\Pi_{12,{m}}:=Q\otimes R^{m},
$$
for the joint distribution of the parameter and the sample, i.e.
$$\Pi_{12,{m}}(A'_{12,m}\times T)=\int_TR_\theta^{m}(A'_{12,m})dQ(\theta),\quad A'_{12,m}\in(\mathcal A_1\times\mathcal  A_2)^{m}, T\in\mathcal T.
$$
The corresponding prior predictive distribution $\beta_{12,{m},Q}^*$ on $(\Omega_1\times\Omega_2)^{m}$ is
$$\beta_{12,{m},Q}^*(A'_{12,m})=\int_\Theta R_\theta^{m}(A'_{12,m})dQ(\theta),\quad A'_{12,m}\in(\mathcal A_1\times\mathcal  A_2)^{m}.
$$
The posterior distribution is a Markov kernel $$R_{m}^*:((\Omega_1\times\Omega_2)^{m},(\mathcal A_1\times\mathcal A_2)^{m})\pt (\Theta,\mathcal T)
$$
such that, for all $A'_{12,m}\in(\mathcal A_1\times\mathcal  A_2)^{m}$ and $T\in\mathcal T$,
$$\Pi_{12,{m}}(A'_{12,m}\times T)=\int_TR_\theta^{m}(A'_{12,m})dQ(\theta)
=\int_{A'_{12,m}}R_{m}^*(x',T)d\beta_{12,{m},Q}^*(x').
$$
Let us write $R_{{m},x'}^*(T):=R_{m}^*(x',T)$.

The posterior predictive distribution on $\mathcal A_1\times\mathcal A_2$  is the Markov kernel
$${R_{{m}}^*}^{R}:((\Omega_1\times\Omega_2)^{m},(\mathcal A_1\times\mathcal  A_2)^m)\pt(\Omega_1\times\Omega_2,\mathcal A_1\times\mathcal  A_2)
$$ 
defined, for $x'\in(\Omega_1\times\Omega_2)^{m}$, by
$${R_{{m}}^*}^{R}(x',A_{12}):=
\int_\Theta R_\theta(A_{12})dR_{{m},x'}^*(\theta)
$$
It follows that, with obvious notations, 
$$\int_{\Omega_1\times\Omega_2}f(x)d{R_{{m},x'}^*}^{\!\!\!\!R}(x)=\int_\Theta\int_{\Omega_1\times\Omega_2}f(x)dR_\theta(x)dR_{{m},x'}^*(\theta)
$$
for any non-negative or integrable real random variable (r.r.v. for short) $f$. 

We can also consider the posterior predictive distribution on $(\mathcal A_1\times\mathcal A_2)^{m}$  defined as the Markov kernel
$${R_{{m}}^*}^{R^{m}}:((\Omega_1\times\Omega_2)^{m},(\mathcal A_1\times\mathcal  A_2)^{m})\pt((\Omega_1\times\Omega_2)^{m},(\mathcal A_1\times\mathcal  A_2)^{m})
$$ 
such that 
$${R_{{m}}^*}^{R^{m}}(x',A'_{12,m}):=
\int_\Theta R_\theta^{m}(A'_{12,m})dR_{{m},x'}^*(\theta)
$$

We introduce some notations for $(x',x,\theta)\in(\Omega_1\times\Omega_2)^{m}\times(\Omega_1\times\Omega_2)\times\Theta$:
\begin{gather*}
\pi'_{m}(x',x,\theta):=x',\quad \pi_{m}(x',x,\theta):=x,\quad \pi_{j,m}(x',x,\theta):=x_j,\;\; j=1,2,\quad
q_{m}(x',x,\theta):=\theta\\ \pi'_{i,m}(x',x,\theta):=x'_i:=(x'_{i1},x'_{i2}),\quad \pi'_{(i),m}(x',x,\theta):=(x'_1,\dots,x'_i),
\end{gather*}
for $1\le i\le m$ (read $i\in\mathbb N$  if  $m=\mathbb N$). 

Let us consider the probability space  
\begin{gather*}
((\Omega_1\times\Omega_2)^m\times(\Omega_1\times\Omega_2)\times\Theta,(\mathcal A_1\times\mathcal A_2)^m\times(\mathcal A_1\times\mathcal A_2)\times\mathcal T,\Pi_m),\qquad \hbox{(2)}
\end{gather*}
where
$$\Pi_m(A'_{12,m}\times A_{12}\times T)=\int_T R_\theta(A_{12}) R_\theta^m(A'_{12,m})dQ(\theta),
$$
when $A'_{12,m}\in(\mathcal A_1\times\mathcal A_2)^m$, $A_{12}\in \mathcal A_1\times\mathcal A_2$ and $T\in\mathcal T$. 

So, for a r.r.v. $f$ on $((\Omega_1\times\Omega_2)^m\times(\Omega_1\times\Omega_2)\times\Theta,(\mathcal A_1\times\mathcal A_2)^m\times(\mathcal A_1\times\mathcal A_2)\times\mathcal T)$,
$$\int f d\Pi_m=\int_\Theta\int_{(\Omega_1\times\Omega_2)^{m}}
\int_{\Omega_1\times\Omega_2} f(x',x,\theta)dR_\theta(x)dR_\theta^m(x')dQ(\theta)\qquad {\rm (3)}
$$
provided that the integral exists. 
Moreover, for a r.r.v. $h$ on $((\Omega_1\times\Omega_2)\times\Theta,(\mathcal A_1\times\mathcal A_2)\times\mathcal T)$,
$$\int h d\Pi_m=\int_\Theta\int_{\Omega_1\times\Omega_2}h(x,\theta)dR_\theta(x)dQ(\theta)
=\int_{\Omega_1\times\Omega_2}\int_\Theta h(x,\theta)dR^*_{1,x}(\theta)d\beta^*_{12,1,Q}(x).
$$ 

The following result is taken from Nogales (2022b).

\begin{myprop}\label{prop1}\rm 
	For $n\in\mathbb N$,
	\begin{gather*}
\Pi_{\mathbb N}^{(\pi'_{(n),\mathbb N},\pi_{1,\mathbb N},q_{\mathbb N})}=\Pi_n,\qquad
\Pi_{\mathbb N}^{(\pi'_{(n),\mathbb N},\pi_{1,\mathbb N})}=\Pi_n^{(\pi'_{(n),n},\pi_{1,n})},\\
\Pi_m^{q_m}=Q,\quad \Pi_m^{(\pi'_m,q_m)}=\Pi_{12,m},\quad \Pi_m^{\pi'_m}=\beta_{12,m,Q}^*,\quad 
\Pi_m^{(\pi_m,q_m)}=\Pi_{12,1},\quad
\Pi_m^{\pi_m}=\beta_{12,1,Q}^*\\ \Pi_m^{\pi'_m|q_m=\theta}=R_\theta^m,\quad \Pi_m^{\pi_m|q_m=\theta}=R_\theta,\quad \Pi_m^{q_m|\pi'_m=x'}=R^*_{m,x'},\quad \Pi_m^{q_m|\pi_m=x}=R^*_{1,x'}
\end{gather*}\end{myprop}

In particular, the probability space (2) contains all the basic ingredients of the Bayesian experiment (1), i.e., the prior distribution, the sampling probabilities, the posterior distributions and the prior predictive distribution. When $m=\mathbb N$, (2) becomes the natural framework to address the asymptotic problem considered in this paper.

\section{Consistency of the Bayes estimator of the regression curve}

Now suppose $(\Omega_2,\mathcal A_2)=(\mathbb R,\mathcal R)$. Let $X_2$ be an squared-integrable r.r.v. such that $E_\theta(X_2^2)$ has a finite prior mean; in particular, $E_\theta(X_2)$ also has a finite prior mean. 

The regression curve of $X_2$ given $X_1$ is the map $x_1\in\Omega_1\mapsto r_\theta(x_1):=E_\theta(X_2|X_1=x_1)$. An estimator of the regression curve $r_\theta$ from a sample of size $n$ of the joint distribution of $(X_1,X_2)$ is a statistic
$$m:(x',x_1)\in(\Omega_1\times \mathbb R)^n\times\Omega_1\longmapsto m(x',x_1)\in\mathbb R,
$$
so that, being observed $x'\in (\Omega_1\times \mathbb R)^n$, $m(x',\cdot)$ is the estimation of $r_\theta$.

From a classical point of view, the simplest way to evaluate the error in estimating an unknown regression curve is to use the expectation of the quadratic deviation (see Nadaraya (1989, p. 120)):
\begin{gather*}
	E_\theta\left[\int_{\Omega_1}(m(x',x_1)-r_\theta(x_1))^2dP_\theta^{X_1}(x_1)\right]=\\
	\int_{(\Omega_1\times \mathbb R)^n}\int_{\Omega_1}(m(x',x_1)-r_\theta(x_1))^2dR_\theta^{p_1}(x_1) dR_\theta^n(x').
\end{gather*}

From a Bayesian point of view, the Bayes estimator of the regression curve $r_\theta$ should minimize the Bayes risk (i.e. the prior mean of the expectation of the quadratic deviation)
\begin{gather*}
	\int_\Theta\int_{(\Omega_1\times \mathbb R)^n}\int_{\Omega_1}(m(x',x_1)-r_\theta(x_1))^2dR_\theta^{p_1}(x_1) dR_\theta^n(x')dQ(\theta)=\\
	E_{\Pi_n}\big[(m(x',x_1)-r_\theta(x_1))^2\big].
\end{gather*}

Recall from Nogales (2022) that the regression curve of $p_2$ on $p_1$ with respect to the posterior predictive distribution ${R^*_{n,x'}}^{\!\!\!\!\!R}$ $$m_n^*(x',x_1):=E_{{R^*_{n,x'}}^{\!\!\!\!\!R}}(p_2|p_1=x_1)
$$ 
is the Bayes estimator of the regression curve $r_\theta(x_1):=E_\theta(X_2|X_1=x_1)$ for the squared error loss function, i.e., 
$$E_{\Pi_n}[(m_n^*(x',x_1)-r_\theta(x_1))^2]\le 
	E_{\Pi_n}[(m_n(x',x_1)-r_\theta(x_1))^2]
	$$
	for any other estimator $m_n$ of the regression curve $r_\theta$.	

We wonder about the consistency of this Bayes estimator. Another question of interest is whether the Bayes risk $E_{\Pi_n}[(m_n^*(x',x_1)-r_\theta(x_1))^2]$ converges to 0 when $n$ goes to $\infty$.

The following result is key to solving the problem. 

\begin{mylemma}\label{lemma2}\rm Let $Y(x',x,\theta):=E_\theta(p_2|p_1=x_1)$.
	
	(i)	For $n\in\mathbb N$ and  we have that
	\begin{gather*}
			E_{{R^*_{n,x'_{(n)}}}^{\!\!\!\!\!\!\!\!\!\!\!\!\!R}}(p_2|p_1=x_1)=E_{\Pi_{\mathbb N}}(Y|(\pi'_{(n),\mathbb N},\pi_{1,\mathbb N})=(x'_{(n)},x_1))
\end{gather*}

	(ii) $$E_{{R^*_{\mathbb N,x'}}^{\!\!\!\!\!\!\!R}}(p_2|p_1=x_1)=E_{\Pi_{\mathbb N}}(Y|(\pi'_{\mathbb N},\pi_{1,\mathbb N})=(x',x_1))
	$$
\end{mylemma}

\begin{myproo}\rm 
	(i)	According to Lemma 1 of Nogales (2022), we have that, for all $A'_{12,n}\in(\mathcal A_1\times\mathcal A_2)^n$ and all $A_i\in\mathcal A_i$, $i=1,2$, 
	\begin{gather*}
		\int_{A'_{12,n}\times A_1\times\Omega_2\times \Theta}R_\theta^{p_2|p_1=x_1}(A_2)d\Pi_n(x',x,\theta)=\\
		\int_{A'_{12,n}\times A_1}\left[{R^{*}_{n,x'}}^{\!\!\!\!R}\right]^{p_2|p_1=x_1}\!\!(A_2)d{\Pi_n}^{(\pi',p_1)}(x',x_1).\qquad (14)
	\end{gather*}
The proof of (i) follows in a standard way from this and Proposition 1 as
$$E_\theta(p_2|p_1=x_1)=\int_{\mathbb R}x_2dR_\theta^{p_2|p_1=x_1}(x_2)\quad\hbox{and}\quad
E_{{R^*_{n,x'_{(n)}}}^{\!\!\!\!\!\!\!\!\!\!R}}\,(p_2|p_1=x_1)=\int_{\mathbb R}x_2d{R^*_{n,x'_{(n)}}}^{\!\!\!\!\!\!\!\!\!R}\,(x_2).
$$

(ii) The proof of (ii) follows from (6) as (i) derives from (14). $\Box$
\end{myproo}

When $\mathcal A'_{(n)}:=(\pi'_{(n),\mathbb N},\pi_{1,\mathbb N})^{-1}((\mathcal A_1\times\mathcal A_2)^n\times\mathcal A_1)$, 
we have that  $(\mathcal A'_{(n)})_n$ is an increasing sequence of  sub-$\sigma$-fields of $(\mathcal A_1\times\mathcal A_2)^{\mathbb N}\times \mathcal A_1$ such that  $(\mathcal A_1\times\mathcal A_2)^{\mathbb N}\times \mathcal A_1=\sigma(\cup_n\mathcal A'_{(n)})$. According to the martingale convergence theorem of Lévy, if $Y$ es $(\mathcal A_1\times\mathcal A_2)^{\mathbb N}\times \mathcal A_1\times\mathcal T$-measurable and $\Pi_{\mathbb N}$-integrable, then 
$$E_{\Pi_{\mathbb N}}(Y|\mathcal A'_{(n)})
$$
converges $\Pi_{\mathbb N}$-a.e. and in $L^1(\Pi_{\mathbb N})$ to $E_{\Pi_{\mathbb N}}(Y|(\mathcal A_1\times\mathcal A_2)^{\mathbb N}\times\mathcal A_1)$. 

Let us consider the measurable function 
$$Y(x',x,\theta):=E_\theta(X_2|X_1=x_1). 
$$
Notice that $E_{\Pi_{\mathbb N}}(Y)=\int_\Theta E_\theta(X_2)dQ(\theta)$, so $Y$ is $\Pi_{\mathbb N}$-integrable. 
Hence, it follows from the aforementioned theorem of Lévy that
$$\lim_n E_{{R^*_{n,x'_{(n)}}}^{\!\!\!\!\!\!\!\!\!\!R}}\,(p_2|p_1=x_1) =
E_{{R^*_{\mathbb N,x'}}^{\!\!\!\!\!\!\!R}}(p_2|p_1=x_1),\quad \Pi_{\mathbb N}-\hbox{a.e.} \qquad (15)
$$
and
$$\lim_n\int_{(\Omega_1\times\Omega_2)^{\mathbb N}\times(\Omega_1\times\Omega_2)\times\Theta}\left|E_{{R^*_{n,x'_{(n)}}}^{\!\!\!\!\!\!\!\!\!\!R}}\,(p_2|p_1=x_1)-E_{{R^*_{\mathbb N,x'}}^{\!\!\!\!\!\!\!R}}(p_2|p_1=x_1)\right| d\Pi_{\mathbb N}(x',x,\theta)=0.\qquad (16)
$$
In the next we will assume the following additional regularity conditions: 
\begin{itemize}
	\item[(i)] $(\Omega_1,\mathcal A_1)$ is a  standard Borel space, 
	
	\item[(ii)] $\Theta$ is a Borel subset of a Polish space and  $\mathcal T$ is its Borel $\sigma$-field, 	and
	
	\item[(iii)] $\{R_\theta\colon \theta\in\Theta\}$ is identifiable.
\end{itemize}

As a consequence of a known theorem of Doob (see Theorem 6.9 and Proposition 6.10 from Ghosal et al. (2017, p. 129, 130)) we have that, for every $x_1\in\Omega_1$, 
$$\lim_n \int_{\Theta}E_{\theta'}(X_2|X_1=x_1)
d\Pi_{\mathbb N}^{q_{\mathbb N}|(\pi'_{(n),\mathbb N},\pi_{1,\mathbb N})=(x'_{(n)},x_1)}(\theta')=E_{\theta}(X_2|X_1=x_1),\quad R_\theta^{\mathbb N}-\hbox{a.e.}
$$
for $Q$-almost every $\theta$. Hence, according to Lemma \ref{lemma2} (i), 
$$\lim_n E_{{R^*_{n,x'_{(n)}}}^{\!\!\!\!\!\!\!\!\!\!R}}\,(p_2|p_1=x_1)=E_{\theta}(X_2|X_1=x_1),\quad R_\theta^{\mathbb N}-\hbox{a.e.}
$$
for $Q$-almost every $\theta$.

In particular,
$$\lim_n E_{{R^*_{n,x'_{(n)}}}^{\!\!\!\!\!\!\!\!\!\!R}}\,(p_2|p_1=x_1)=E_{\theta}(X_2|X_1=x_1),\quad \Pi_{\mathbb N}-\hbox{a.e.}
$$
In this sense we can say that the predictive posterior regression curve $E_{{R^*_{n,x'_{(n)}}}^{\!\!\!\!\!\!\!\!\!\!R}}\,(p_2|p_1=x_1)$ of $X_2$ given $X_1=x_1$ is a strongly consistent estimator of the sampling regression curve $E_{\theta}(X_2|X_1=x_1)$ of $X_2$ given $X_1=x_1$. 

From this and (15) we obtain that
$$ E_{{R^*_{\mathbb N,x'}}^{\!\!\!\!\!\!\!R}}\,(p_2|p_1=x_1)=E_{\theta}(X_2|X_1=x_1)
,\quad \Pi_{\mathbb N}-\hbox{a.e.}
$$

According to (16) we obtain that
$$\lim_n\int_{(\Omega_1\times\Omega_2)^{\mathbb N}\times(\Omega_1\times\Omega_2)\times\Theta}\left|E_{{R^*_{n,x'_{(n)}}}^{\!\!\!\!\!\!\!\!\!\!R}}\,(p_2|p_1=x_1)-E_{\theta}(X_2|X_1=x_1)\right| d\Pi_{\mathbb N}(x',x,\theta)=0,\qquad (17)
$$
which proves that the Bayes risk of the optimal estimator $E_{{R^*_{n,x'_{(n)}}}^{\!\!\!\!\!\!\!\!\!\!R}}\,(p_2|p_1=x_1)$ of the regression curve $E_{\theta}(X_2|X_1=x_1)$ converges to 0 for the $L^1$-loss function. 

We wonder if that also happens for the $L^1$-squared loss function, i.e., if the Bayes risk 
$$E_{\Pi_n}[(m_n^*(x',x_1)-r_\theta(x_1))^2]
$$
converges vers 0 when $n$ goes to $\infty$. Theorem 6.6.9 of Ash et al. (2000) shows that the answer is affirmative because 
$$m_n^*(x',x_1)=E_{\Pi_{\mathbb N}}(Y|\mathcal A'_{(n)})
$$
and, by Jensen's inequality,
\begin{gather*}
E_{\Pi_{\mathbb N}}(E_{\Pi_{\mathbb N}}(Y|\mathcal A'_{(n)})^2)\le
E_{\Pi_{\mathbb N}}(E_{\Pi_{\mathbb N}}(Y^2|\mathcal A'_{(n)}))=
E_{\Pi_{\mathbb N}}(Y^2)=\int_\Theta X_2^2dQ(\theta)<\infty.
\end{gather*}

So, we have proved the following result (in fact part (i) is shown in Nogales (2022a) and its statement is reproduced here for the sake of completeness).

\begin{mytheo}\rm Let $(\Omega,\mathcal A,\{P_\theta\colon\theta\in(\Theta,\mathcal T,Q)\})$ be a Bayesian statistical experiment, and $X_1:(\Omega,\mathcal A,\{P_\theta\colon\theta\in(\Theta,\mathcal T,Q)\})\rightarrow(\Omega_1,\mathcal A_1)$ and $X_2:(\Omega,\mathcal A,\{P_\theta\colon\theta\in(\Theta,\mathcal T,Q)\})\rightarrow(\mathbb R,\mathcal R)$ two statistics such that $E_\theta(X_2^2)$ has finite prior mean. 
Let us suppose that conditions (i)-(iii) above hold. 	
 Then:
 
 (i) The regression curve of $p_2$ on $p_1$ with respect to the posterior predictive distribution ${R^*_{n,x'}}^{\!\!\!\!\!R}$ $$m_n^*(x',x_1):=E_{{R^*_{n,x'}}^{\!\!\!\!\!R}}(p_2|p_1=x_1)$$ is the Bayes estimator of the regression curve $r_\theta(x_1):=E_\theta(X_2|X_1=x_1)$ for the squared error loss function, i.e., $$E_{\Pi_n}[(m_n^*(x',x_1)-r_\theta(x_1))^2]\le 
 E_{\Pi_n}[(m_n(x',x_1)-r_\theta(x_1))^2]
 $$
 for any other estimator $m_n$ of the regression curve $r_\theta$.	 
 
 (ii) Moreover, $m_n^*$ is a strongly consistent estimator of thr regression curve, in the sense that
 $$\lim_n E_{{R^*_{n,x'_{(n)}}}^{\!\!\!\!\!\!\!\!\!\!R}}\,(p_2|p_1=x_1)=E_{\theta}(X_2|X_1=x_1),\quad \Pi_{\mathbb N}-\hbox{a.e.}
 $$
 
 (iii) Finally, the Bayes risk of $m_n^*$ converges to 0 both  for the $L^1$-loss function and the $L^1$-squared loss function, i.e.,
 $$\lim_n E_{\Pi_{\mathbb N}}[|m_n^*(x',x_1)-r_\theta(x_1)|^k]=0,\quad k=1,2.
 $$
 \end{mytheo}

\section{Examples.}\label{secex}

\begin{myexa}\rm Let us suppose that, for $\theta,\lambda,x_1>0$, $P_\theta^{X_1}= G(1,\theta^{-1})$, $P_\theta^{X_2|X_1=x_1}=G(1,(\theta x_1)^{-1})$ and $Q=G(1,\lambda^{-1})$, where $G(\alpha,\beta)$ denotes the gamma distribution of parameters $\alpha,\beta>0$. Hence the joint density of $X_1$ and $X_2$ is
	$$f_\theta(x_1,x_2)=\theta^2 x_1\exp\{-\theta x_1(1+x_2)\}I_{]0,\infty[^2}(x_1,x_2).
	$$
	It is shown in Nogales (2022b), Example 1, the Bayes estimator of the regression curve $r_\theta(x_1):=E_\theta(X_2|X_1=x_1)=\frac1{\theta x_1}$ is, for $x_1>0$, 
	$$m^*_n(x',x_1)=\int_0^\infty x_2\cdot{f^*_{n,x'}}^{\!\!\!\!X_2|X_1=x_1}(x_2)dx_2=\frac{\lambda+x_1+\sum_{i=1}^nx'_{i1}(1+x'_{i2})}{(2n+1)x_1}.
	$$
	Theorem 1 shows that this a strongly consistent estimator of the regression curve $r_\theta(x_1)$.
\end{myexa} 

\begin{myexa}\rm Let us suppose that $X_1$ has a Bernoulli distribution of unknown parameter $\theta\in]0,1[$ (i.e. $P_\theta^{X_1}=Bi(1,\theta)$) and, given $X_1=k_1\in\{0,1\}$, $X_2$ has distribution $Bi(1,1-\theta)$ when $k_1=0$ and $Bi(1,\theta)$ when $k_1=1$, i.e. $P_\theta^{X_2|X_1=k_1}=Bi(1,k_1+(1-2k_1)(1-\theta))$. We can think of tossing a coin with probability $\theta$ of getting heads ($=1$) and making a second toss of this coin if it comes up heads on the first toss, or tossing a second coin with probability $1-\theta$ of making heads if the first toss is tails ($=0$). Consider the uniform distribution on $]0,1[$ as the prior distribution $Q$.
	
	So, the joint probability function of $X_1$ and $X_2$ is
	\begin{gather*}\begin{split}
			f_\theta(k_1,k_2)&=\theta^{k_1}(1-\theta)^{1-k_1}[k_1+(1-2k_1)(1-\theta)]^{k_2}[1-k_1-(1-2k_1)(1-\theta)]^{1-k_2}\\
			&=\begin{cases}
				\theta(1-\theta)\hbox{ if } k_2=0,\\ (1-\theta)^2\hbox{ if } k_1=0, k_2=1,\\ \theta^2\hbox{ if } k_1=1, k_2=1.
			\end{cases}
	\end{split}\end{gather*}
	It is shown in Nogales(2022b), Example 2, that 
	the Bayes estimator of the conditional mean $r_\theta(k_1):=E_\theta(X_2|X_1=k_1)=\theta^{k_1}(1-\theta)^{1-k_1}$ is, for $k_1=0,1$,
	$$m^*_n(k',k_1)={f^*_{n,k'}}^{\!\!\!\!X_2|X_1=k_1}(1)=\begin{cases}
		\frac{n_{+0}(k')+2n_{01}(k')+1}{2n+n_{+0}(k')+2n_{01}(k')+3}&\hbox{ if } k_1=0,\vspace{1ex}\\
		\frac{n_{+0}(k')+2n_{01}(k')+1}{2n+n_{+0}(k')+2n_{01}(k')+4}&\hbox{ if } k_1=1,
	\end{cases}	
	$$
	being $n_{j_1j_2}(k')$ the number of indices $i\in\{1,\dots,n\}$ such that $(k'_{i1},k'_{i2})=(j_1,j_2)$ and $n_{+j}=n_{0j}+n_{1j}$ for $j=0,1$. 
	
	Theorem 1 proves that it is a strongly consistent estimator of the conditional mean $r_\theta(k_1)$.
\end{myexa}

\begin{myexa}\rm 
	Let $(X_1,X_2)$ have bivariate normal distribution  
	$$N_2\left(
	\left(\begin{array}{c}
		\theta\\
		\theta
	\end{array}\right),\sigma^2
	\left(\begin{array}{cc}
		1 & \rho\\
		\rho & 1
	\end{array}\right)
	\right),
	$$
	and consider the prior distribution $Q=N(\mu,\tau^2)$. 
	It is shown in Nogales (2022b), Example 3, that
	that the conditional mean
	$$E_{{R^*_{n,x'}}^{\!\!\!\!\!R}}(p_2|p_1=x_1)=(1-\rho_1)m_1(x')+\rho_1x_1
	$$
	is the Bayes estimator of the regression curve 
	$$E_\theta(X_2|X_1=x_1)=(1-\rho)\theta+\rho x_1
	$$
	for the squared error loss function, where
	\begin{gather*}	
		\rho_1=-\frac{a_n(\rho,\sigma,\tau)+\frac{1-\rho}{1+\rho}}{a_n(\rho,\sigma,\tau)-\frac{1-\rho}{1+\rho}}\cdot\rho,
		\quad 
		m_1(x')=\frac{s_1(x')+(1+\rho)\frac{\sigma^2}{\tau^2}\mu}{2(1-\rho_1)(1+\rho)^2\sigma^2a_n(\rho,\sigma,\tau)},
	\end{gather*}	
	being
	\begin{gather*}s_1(x'):=\sum_i(x'_{i1}+x'_{i2}),\quad
		a_n(\rho,\sigma,\tau):=2(n+1)(1+\rho)+\frac{\sigma^2}{\tau^2}.
		\end{gather*}	
		
	Theorem 1 proves that it is a strongly consistent estimator of this regression curve. 
\end{myexa}

\section{Appendix.}

We recover here the Appendix of Nogales (2022a) to briefly recall some basic concepts about Markov kernels, mainly to fix the notations. In the next,  $(\Omega,\mathcal A)$, $(\Omega_1,\mathcal A_1)$ and so on will denote measurable spaces. 

\begin{mydef}\rm  1) (Markov kernel) A Markov kernel
	$M_1:(\Omega,\mathcal A)\pt   (\Omega_1,\mathcal A_1)$ is a map $M_1:\Omega\times\mathcal A_1\rightarrow[0,1]$ such that: 	(i) $\forall \omega\in\Omega$, $M_1(\omega,\cdot)$ is a  probability
	measure on 	$\mathcal A_1$, (ii) $\forall A_1\in\mathcal A_1$, $M_1(\cdot,A_1)$ is $\mathcal A$-measurable.\par
	2)	(Image of a Markov kernel) The image (or {\it probability
		distribution}) of a Markov kernel $M_1:(\Omega,\mathcal A,P)\pt
	(\Omega_1,\mathcal A_1)$ on a probability space is the probability
	measure  $P^{M_1}$ on $\mathcal A_1$ defined by
	$P^{M_1}(A_1):=\int_{\Omega}M_1(\omega,A_1)\,dP(\omega)$.
	\par
	3)  (Composition of Markov kernels) Given two Markov kernels
	$M_1:(\Omega_1,\mathcal A_1)\pt (\Omega_2,\mathcal A_2)$ and
	$M_2:(\Omega_2,\mathcal A_2)\pt (\Omega_3,\mathcal A_3)$, its composition  is defined 	as the Markov kernel $M_2M_1:(\Omega_1,\mathcal A_1)\pt
	(\Omega_3,\mathcal A_3)$ given by
	$$M_2M_1(\omega_1,A_3)=\int_{\Omega_2}M_2(\omega_2,A_3)M_1(\omega_1,d\omega_2).
	$$
\end{mydef}

\begin{myprem}\rm  1) (Markov kernels as extensions of the concept of random variable) The concept of Markov kernel extends the concept of random variable (or measurable map). A random variable $T_1:(\Omega,\mathcal 	A,P)\rightarrow(\Omega_1,\mathcal A_1)$ will be identified with the Markov kernel $M_{T_1}:(\Omega,\mathcal A;P)\pt   (\Omega_1,\mathcal
	A_1)$ defined by $M_{T_1}(\omega,A_1)=\delta_{T_1(\omega)}(A_1)=I_{A_1}(T_1(\omega))$,
	where $\delta_{T_1(\omega)}$ denotes the Dirac measure -the
	degenerate distribution- at the point $T_1(\omega)$, and $I_{A_1}$ is
	the indicator function of the event $A_1$. In particular, the probability distribution $P^{M_{T_1}}$ of $M_{T_1}$ coincides with the probability distribution $P^{T_1}$ of $T_1$ defined as $P^{T_1}(A_1):=P(T_1\in A_1)$\par	
	2) Given a Markov kernel $M_1:(\Omega_1,\mathcal A_1)\pt (\Omega_{2},\mathcal A_{2})$ and a random variable $X_2:(\Omega_2,\mathcal A_2)\rightarrow (\Omega_{3},\mathcal A_{3})$, we have that  $M_{X_2}M_1(\omega_1,A_3)=M_1(\omega_1,X_2^{-1}(A_3))=
	M_1(\omega_1,\cdot)^{X_2}(A_3).$ We write $X_2M_1:=M_{X_2}M_1$.\par
	3) Given a Markov kernel $M_1:(\Omega_1,\mathcal A_1,P_1)\pt (\Omega_{2},\mathcal A_{2})$ we write $P_1\otimes M_1$ for the only probability measure on the product $\sigma$-field $\mathcal A_1\times\mathcal A_2$ such that 
	$$(P_1\otimes M_1)(A_1\times A_2)=\int_{A_1}M_1(\omega_1,A_2)dP_1(\omega_1),\quad A_i\in\mathcal A_i,\, i=1,2.$$
	
	4) Given two r.v. $X_i:(\Omega,\mathcal A,P)\rightarrow(\Omega_i,\mathcal A_i)$, $i=1,2$, we write $P^{X_2|X_1}$ for the conditional distribution of $X_2$ given $X_1$, i.e. for the Markov kernel $P^{X_2|X_1}:(\Omega_1,\mathcal A_1)\pt(\Omega_2,\mathcal A_2)$ such that $$P^{(X_1,X_2)}(A_1\times A_2)=\int_{A_1}P^{X_2|X_1=x_1}(A_2)dP^{X_1}(x_1),\quad  A_i\in\mathcal A_i,\,i=1,2.	
	$$
	So $P^{(X_1,X_2)}=P^{X_1}\otimes P^{X_2|X_1}$. $\Box$
	
\end{myprem}

Let $(\Omega,\mathcal A,\{P_\theta\colon\theta\in(\Theta,\mathcal T,Q)\})$ be a Bayesian statistical experiment, 
where $Q$ denotes the prior distribution on the parameter space $(\Theta,\mathcal T)$. We suppose that $P(\theta,A):=P_\theta(A)$ is a Markov kernel $P:(\Theta,\mathcal T)\pt(\Omega,\mathcal A)$. When needed we shall suppose that $P_\theta$ has a density (Radon-Nikodym derivative) $p_\theta$ with respect to a $\sigma$-finite measure $\mu$ on $\mathcal A$ and that the likelihood function $\mathcal L(\omega,\theta):=p_\theta(\omega)$ is $\mathcal A\times\mathcal T$-measurable (this is sufficient to prove that $P$ is a Markov kernel). 

Let $\Pi:=Q\otimes P$, i.e.
$$\Pi(A\times T)=\int_TP_\theta(A)dQ(\theta), \quad A\in\mathcal A, T\in\mathcal T.
$$
The prior predictive distribution is $\beta_Q^*:=\Pi^I$ (the distribution of $I$ with respect to $\Pi$), where $I(\omega,\theta):=\omega$. So
$$\beta_Q^*(A)=\int_\Theta P_\theta(A)dQ(\theta).
$$
The posterior distribution is a Markov kernel $P^*:(\Omega,\mathcal A)\pt(\Theta,\mathcal T)$ such that
$$\Pi(A\times T)=\int_TP_\theta(A)dQ(\theta)=\int_AP^*_\omega(T)d\beta_Q^*(\omega), \quad A\in\mathcal A, T\in\mathcal T,
$$
i.e. such that $\Pi=Q\otimes P=\beta_Q^*\otimes P^*$. This way the Bayesian statistical experiment can be identified with the probability space $(\Omega\times\Theta,\mathcal A\times\mathcal T,\Pi)$, as proposed, for instance, in Florens et al. (1990). 

It is well known that, for $\omega\in\Omega$, the posterior $Q$-density is proportional to the likelihood, i.e.
$$p^*_{\omega}(\theta):=\frac{dP^*_\omega}{dQ}(\theta)=C(\omega)p_\theta(\omega)
$$
where $C(\omega)=[\int_\Theta p_\theta(\omega)dQ(\theta)]^{-1}$.

The posterior predictive distribution on $\mathcal A$ given $\omega$ is 
$${P_\omega^*}^P(A)=\int_\Theta P_\theta(A)dP_\omega^*(\theta),\quad A\in\mathcal A.
$$
This is a Markov kernel 
$$PP^*(\omega,A):={P_\omega^*}^P(A).$$
It is readily shown that the posterior predictive density is
$$\frac{d{P_\omega^*}^P}{d\mu}(\omega')=\int_\Theta p_\theta(\omega')p^*_\omega(\theta)dQ(\theta).
$$

We know from Nogales (2022a) that
$$\int_{\Omega\times\Theta}\sup_{A\in\mathcal A}|{P_\omega^*}^P(A)-P_\theta(A)|^2d\Pi(\omega,\theta)\le
\int_{\Omega\times\Theta}\sup_{A\in\mathcal A}|M(\omega,A)-P_\theta(A)|^2d\Pi(\omega,\theta),$$
for every Markov kernel $M:(\Omega,\mathcal A)\pt (\Omega,\mathcal A)$ provided that $\mathcal A$ is separable  (recall that a $\sigma$-field is said to be separable, or countably generated, if it contains a countable subfamily which generates it). We also have that, for a real statistic  $X$  with finite mean, the posterior predictive mean
$$E_{(P_\omega^*)^P}(X)=\int_\Theta\int_\Omega X(\omega')dP_\theta(\omega')dP_\omega^*(\theta)$$
is the Bayes estimator  of $f(\theta):=E_\theta(X)$, as $E_{(P_\omega^*)^P}(X)=E_{P_\omega^*}(E_\theta(X))$.

\section{Acknowledgements.}
This paper has been supported by the Junta de Extremadura (Spain) under the grant Gr21044.
\vspace{1ex}

\section{References.}

\begin{itemize}

	\item Ash, R.B., Dóleans-Dade, C. (2000) Probability and Measure Theory, 2nd Ed., Academic Press, San Diego USA.

	\item Geisser, S. (1993) Predictive Inference: An Introduction, Springer Science+ Business Media, Dordrecht.

	\item Gelman, A., Carlin, J.B., Stern, H.S., Dunson, D.B., Vehtari, A., Rubin, D.B. (2014) Bayesian Data Analysis, 3rd ed., CRC Press.
	
	\item Ghosal, S., Vaart, A.v.d. (2017) Fundamentals of Noparametric Bayesian Inference, Cambridge University Press, Cambridge UK. 
	
	  \item Nadaraya, E.A. (1989), Nonparametric Estimation of Probability Densities and Regression Curves, Kluwer Academic Publishers, Dordrecht. 	
	
	\item Nogales, A.G. (2022a), On Bayesian Estimation of Densities and Sampling Distributions: the Posterior Predictive Distribution as the Bayes Estimator, Statistica Neerlandica 76(2), 236-250.
	
	\item Nogales, A.G. (2022b), Optimal Bayesian Estimation of a Regression Curve, a Conditional Density, and a Conditional Distribution, Mathematics 10(8), 1213.
	
	\item Nogales, A.G. (2022c), On consistency of the Bayes Estimator of the Density, Mathematics, 10(4), 636.
	
\end{itemize}

\end{document}